\documentclass[12pt]{amsart}

\usepackage{amsmath,amssymb,amscd,amsxtra,upref,stmaryrd}
\usepackage{latexsym}
\usepackage[dvips]{graphics,epsfig}
\usepackage{enumerate,mathrsfs}
\usepackage[colorlinks=true, pdfstartview=FitV, linkcolor=blue, citecolor=blue, urlcolor=blue]{hyperref}
\usepackage{url}

\allowdisplaybreaks

\numberwithin{equation}{section}

\newtheorem{theorem}{Theorem}[section]

\newtheorem{proposition}[theorem]{Proposition}
\newtheorem{corollary}[theorem]{Corollary}

\newtheorem{remark}[theorem]{Remark}

\newcommand{\supp}{\text{supp}}

\newcommand{\beq}{\begin{equation}}
\newcommand{\eeq}{\end{equation}}

\newcommand{\bal}{\begin{align*}}
\newcommand{\eal}{\end{align*}}

\newcommand{\subR}{{\mathbb R}}
\newcommand{\subRn}{{{\mathbb R}^n}}

\newcommand{\nR}{{\mathbb R}}
\newcommand{\loc}{{\mathrm loc}}

\DeclareMathOperator*{\essinf}{ess\,inf}
\DeclareMathOperator*{\esssup}{ess\,sup}
\DeclareMathOperator{\sgn}{sgn}

\newcommand{\pp}{{p(\cdot)}}
\newcommand{\cpp}{{p'(\cdot)}}

\newcommand{\qq}{{q(\cdot)}}
\newcommand{\cqq}{{q'(\cdot)}}

\def\Xint#1{\mathchoice
   {\XXint\displaystyle\textstyle{#1}}%
   {\XXint\textstyle\scriptstyle{#1}}%
   {\XXint\scriptstyle\scriptscriptstyle{#1}}%
   {\XXint\scriptscriptstyle\scriptscriptstyle{#1}}%
   \!\int}
\def\XXint#1#2#3{{\setbox0=\hbox{$#1{#2#3}{\int}$}
     \vcenter{\hbox{$#2#3$}}\kern-.5\wd0}}

\def\avgint{\Xint-}

\begin{document}

\title[commutators on Banach function spaces]{Necessary conditions for the boundedness of linear and bilinear commutators on Banach function spaces}

\author{Lucas Chaffee}
\address{Department of Mathematics,
Western Washington University, Bellingham, WA 98225-9063}
\email{lucas.chaffee@wwu.edu}

\author{David Cruz-Uribe, OFS}
\address{Department of Mathematics\\
University of Alabama, Box 870350, Tuscaloosa, AL 35487}
\email{dcruzuribe@ua.edu}




\date{January 26, 2017}

\thanks{ The second
  author is supported by NSF Grant DMS-1362425 and research funds from the
  Dean of the College of Arts \& Sciences, the University of Alabama.}

\begin{abstract}
  In this article we extend recent results by the first
  author~\cite{LC} on the necessity of $BMO$ for the boundedness of
  commutators on the classical Lebesgue spaces.  We generalize these
  results to a large class of Banach function spaces.  We show that
  with modest assumptions on the underlying spaces and on the operator
  $T$, if the commutator $[b,T]$ is bounded, then the function $b$ is
  in $BMO$.

\end{abstract}

\maketitle

\section{Introduction}

The purpose of this paper is to extend a recent result of the first
author~\cite{LC} on necessary conditions for commutators to be bounded
on the classical Lebesgue spaces.  He showed that if $T$ is a ``nice''
operator, and if (for example) the commutator $[b,T]$ is bounded on
$L^p$, then $b\in BMO$.  He also proved an analogous result for bilinear
commutators.  We generalize these results to a large collection of
Banach function spaces.  To do so requires the assumption of a
geometric condition on the underlying spaces that is closely related
to the boundedness of the Hardy-Littlewood maximal operator, and which
holds in a large number of important special cases.

To state our results we recall some basic facts about Banach function
spaces.  For further information, see Bennett and
Sharpley~\cite{bennett-sharpley88}.  By a Banach function space $X$ we
mean a Banach space of measurable functions over $\nR^n$ whose norm
$\|\cdot\|_X$ satisfies the following for all $f,\,g\in X$:

\begin{enumerate}
\setlength{\itemsep}{3pt}

\item $\|f\|_X = \||f|\|_X$;
\item  if  $|f|\le |g|$ a.e., then
  $\|f\|_X\leq \|g\|_X$;

\item  if $\{f_n\}\subset X$ is a
  sequence such that $|f_n|$ increases to $|f|$ a.e., then
  $\|f_n\|_X$ increases to $\|f\|_X$;

\item if $E\subset \nR^n$ is bounded, then $\|\chi_E\|_X<\infty$;

\item if $E$ is bounded, then $\int_E |f(x)|\, d\mu\le C\|f\|_X$,
  where $C=C(E,X)$.
\end{enumerate}

Given a Banach function space $X$, there exists another Banach
function space $X'$, called the associate space of $X$, such that  for
all $f\in X$,
\[ \|f\|_X \approx \sup_{\substack{g\in X'\\ \|g\|_{X'}\leq 1}}
\int_\subRn f(x)g(x)\,dx. \]
In many (but not all) cases, the associate space is equal to the dual
space $X^*$.  The associate space, however, is always reflexive, in
the sense that $(X')'=X$.   Moreover, we have the following
generalization of H\"older's inequality:
\[ \int_{\nR^n} |f(x)g(x)|\,dx \lesssim \|f\|_X \|g\|_{X'}. \]

\medskip

Given a linear operator $T$,
define the commutator $[b,T]f(x) = b(x)Tf(x)-T(bf)(x)$, where $b$ is a
locally integrable function.  We can now state our first result.

\begin{theorem} \label{thm:linear} Given Banach function spaces $X$
  and $Y$, and $0\leq \alpha<n$, suppose that for every cube $Q$,
\begin{equation} \label{eqn:linear1}
 |Q|^{-\frac{\alpha}{n}}\|\chi_Q\|_{Y'} 
\|\chi_Q\|_{X} \leq C|Q|. 
\end{equation}
 Let $T$ be a linear operator defined on $X$ which can
 be represented by  
\[ Tf(x)=\int_\subRn K(x-y)f(y)\,dy\]
for all $x\not\in  \supp(f)$, where $K$ is a homogeneous
kernel of degree $n-\alpha$.  Suppose further that there exists a
ball $B\in\subR^{n}$ on which $\frac1K$ has an absolutely convergent
Fourier series.   If the commutator satisfies
$[b,T] : X\to Y$, 
then $b \in BMO(\nR^n)$.   
\end{theorem}

A wide variety of classical operators satisfy the hypotheses of
Theorem~\ref{thm:linear}.  The kernel $K$ is such that $\frac1K$ has
an absolutely convergent Fourier series if it is non-zero on a ball
$B$ and has enough regularity: $K\in C^s(B)$ for $s>n/2$ is
sufficient.  (See Grafakos~\cite [Theorem~3.2.16]{grafakos08a}.  For
weaker sufficient conditions, see recent results by M\'oricz and
Veres~\cite{MR2361606}.)  In the linear case this condition is
satisfied by Calder\'on-Zygmund singular integrals of convolution type
whose kernels are smooth, and in particular the Riesz transforms.  It
also includes the fractional integral operators (also referred to as
Riesz potentials).  For precise definitions, see
Section~\ref{section:bfs} below.

\medskip

To state our  result in bilinear case we recall that there are two
commutators to consider:  if $T$ is a bilinear operator
and  $b\in L^1_\loc(\nR^n)$, define
\begin{gather*}
 [b,T]_1(f,g)(x) =  b(x)T(f,g)(x)- T(bf,g)(x), \\
[b,T]_2(f,g)(x) =  b(x)T(f,g)(x)- T(f,bg)(x).
\end{gather*}

\begin{theorem} \label{thm:main}
Given  Banach function spaces $X_1,\ X_2$, and $Y$, and $0\leq
\alpha<2n$, suppose  that for every cube $Q$,
\begin{equation} \label{eqn:main1}
 |Q|^{-\frac{\alpha}{n}}\|\chi_Q\|_{Y'} 
\|\chi_Q\|_{X_1}\|\chi_Q\|_{X_2} \leq C|Q|. 
\end{equation}
 Let $T$ be a bilinear operator defined on $X_1\times X_2$ which can
 be represented by  
\[ T(f,g)(x)=\int_\subRn\int_\subRn K(x-y,z-y)f(y)g(z)\,dydz \]
for all $x\not\in  \supp(f)\cap \supp(g)$, where $K$ is a homogeneous
kernel of degree $2n-\alpha$.  Suppose further that there exists a
ball $B\in\subR^{2n}$ on which $\frac1K$ has an absolutely convergent
Fourier series.   If for $j=1$ or $j=2$, the bilinear commutator satisfies
 \[   [b,T]_j:X_1\times X_2\to Y, \]
then $b \in BMO(\nR^n)$.   
\end{theorem}

\begin{remark}
Theorem~\ref{thm:main} extends naturally to 
multilinear operators.  We leave the statement and proof of this
generalization to the interested reader.
\end{remark}

\begin{remark}
The restrictions on $\alpha$ are not actually necessary in the proofs
of Theorems~\ref{thm:linear} and~\ref{thm:main}:
we can take any $\alpha \in \nR$.  However, we are not aware of any
operators for which Banach function space estimates hold that do not
satisfy the given restrictions on $\alpha$.
\end{remark}

\medskip

For the absolute convergence of multiple Fourier series, see the above
references.  In the bilinear case, Theorem~\ref{thm:main} covers such
operators as the bilinear Calder\'on-Zygmund operators with smooth
kernels~\cite{MR2030573,MR1880324,MR1947875,MR2483720} and the
bilinear fractional integral
operator~\cite{MR1164632,MR1812822,MR1682725,MR2514845}.

One drawback of Theorem~\ref{thm:main} is that we must assume that the
target space $Y$ is a Banach function space.  This is somewhat restrictive:
even in the case of the Lebesgue spaces, bilinear operators
satisfy inequalities of the form $T : L^{p_1}\times L^{p_2}
\rightarrow L^p$ where $p<1$.    This assumption, however, is
intrinsic to the statement and proof of our result, since we need to
use the generalized reverse H\"older inequality.   We are uncertain
what the correct assumption should be when we assume that $Y$ is only
a quasi-Banach space.

\begin{remark} 
The necessity of BMO for the boundedness on Lebesgue spaces of
  commutators of certain multilinear singular integrals was recently
  proved in~\cite{LW} using a completely different approach, but the
  authors were also required to assume that for the target space
  $L^p$, $p\geq 1$.  
\end{remark}

\medskip

The remainder of this paper is split into two parts.  We defer the
actual proof of Theorems~\ref{thm:linear} and~\ref{thm:main} to
Section~\ref{section:proof} and in fact we will only give the proof of
the latter; the proof of the former is gotten by a trivial adaptation
of the proof in the bilinear case.  But first, in Section~\ref{section:bfs} we
give a number of specific examples of Banach function spaces and
consider the relationship between the conditions~\eqref{eqn:linear1}
and ~\eqref{eqn:main1}, and sufficient conditions for maximal
operators and  commutators to
be bounded.  

Throughout this paper, $n$ will denote the dimension of the underlying
space, $\nR^n$.  We will consider real-valued functions over
$\nR^n$.  Cubes in $\nR^n$ will always have their sides parallel to
the coordinate axes.   If we write $A\lesssim B$, we mean $A\leq cB$, where the
constant $c$ depends on the operator $T$, the Banach function spaces,
and on the dimension $n$.  These implicit constants may change from
line to line.   If we write $A\approx B$, then $A\lesssim B$ and
$B\lesssim A$.  

\section{Examples of Banach function spaces}
\label{section:bfs}

\subsection*{Averaging and maximal operators}
The necessary condition \eqref{eqn:linear1} in Theorem~\ref{thm:linear}
is closely related to a necessary condition for the boundedness of
averaging operators and fractional maximal operators.  For
$0\leq \alpha<n$, given a cube $Q$ define the linear $\alpha$-averaging operator
\[ A_\alpha^Q f(x) = |Q|^{\frac{\alpha}{n}}\avgint_Q f(y)\,dy
\cdot \chi_Q(x).  \]
We define the associated fractional maximal operator by
\[ M_\alpha f(x) = \sup_Q |Q|^{\frac{\alpha}{n}}\avgint_Q |f(y)|\,dy
\cdot \chi_Q(x). \]
We immediately have that for all $Q$, $|A_\alpha^Qf(x)|\leq M_\alpha
f(x)$.  We also make the analogous definitions in the bilinear case:
for $0\leq \alpha <2n$, 
\[ A_\alpha^Q (f,g)(x) = |Q|^{\frac{\alpha}{n}}\avgint_Q f(y)\,dy
\avgint_Q g(y)\,dy\cdot \chi_Q(x),  \]
\[ M_\alpha (f,g)(x) = \sup_Q |Q|^{\frac{\alpha}{n}}\avgint_Q |f(y)|\,dy
\avgint_Q |g(y)|\,dy\cdot \chi_Q(x).\]
Again, we have the pointwise bound $|A_\alpha^Q(f,g)(x)|\leq M_\alpha
(f,g)(x)$.   In both the linear and bilinear case, when $\alpha=0$ we
write $M$ instead of $M_0$.

In the linear case the maximal operators are classical; the averaging
operators were implicit but seem to have first been considered
explicitly when $\alpha=0$ in~\cite{jawerth86}.   In the bilinear
case, when $\alpha=0$,
the bilinear maximal operator was introduced in~\cite{MR2483720}, and
when $0<\alpha<2n$ in~\cite{MR2514845}.   The bilinear averaging
operators were first considered in~\cite{Kokilashvili:2015gw}.  
The following result is due to Berezhnoi~\cite[Lemma~2.1]{MR1622773}
in the linear case and to Kokilashvili {\em et
  al.}~\cite[Theorem~2.1]{Kokilashvili:2015gw} in the bilinear case.  

\begin{proposition} \label{prop:avg-op}
Fix $0\leq \alpha<n$.  Given Banach function spaces $X$, $Y$, 
there exists a constant $C$ such that for every cube $Q$, 
\begin{equation} \label{eqn:avg-op1}
  \|\chi_Q\|_{Y}\|\chi_Q\|_{X'}  \leq C|Q|^{1-\frac{\alpha}{n}};
\end{equation}
if and only if  
\begin{equation} \label{eqn:avg-op2}
 \|A_\alpha^Q f\|_{Y} \leq C\|f\|_{X}. 
\end{equation}

Similarly, in the bilinear case fix $0\leq \alpha <2n$.  Given Banach
function spaces $X_1$, $X_2$ and $Y$ there exists a constant $C$ such
that for every cube $Q$,
\begin{equation} \label{eqn:bi-avg-op1}
  \|\chi_Q\|_{Y}\|\chi_Q\|_{X_1'}\|\chi_Q\|_{X_2'}  
\leq C|Q|^{1-\frac{\alpha}{n}};
\end{equation}
if and only if 
\begin{equation} \label{eqn:bi-avg-op2}
 \|A_\alpha^Q (f,g)\|_{Y} \leq C\|f\|_{X_1}\|g\|_{X_2}. 
\end{equation}
\end{proposition}

By the pointwise estimates above,~\eqref{eqn:avg-op2} holds whenever
the fractional maximal operator satisfies
$M_\alpha : X \rightarrow Y$, and the corresponding fact is true in
the bilinear case.  Moreover, when $\alpha=0$ and $X=Y$,
condition~\eqref{eqn:avg-op1} is the same as~\eqref{eqn:linear1}.
This yields the following important corollary to Theorem~\ref{thm:linear}.

\begin{corollary} \label{cor:m-bounded}
Given a Banach function space $X$, if $M : X \rightarrow X$,
then~\eqref{eqn:linear1} holds.  Equivalently, if the maximal operator
is bounded and $T$ is an operator with a kernel that is homogenous of
degree $n$, then a necessary condition for $[b,T] : X \rightarrow X$
is that $b\in BMO$.   
\end{corollary}

\medskip

As we will discuss below, the assumption that the maximal operator is
bounded on the Banach function space $X$ is a natural one.
Unfortunately, we cannot generalize Corollary~\ref{cor:m-bounded} to
the case $\alpha>0$ for linear operators or to any bilinear operators
acting on general Banach function spaces.  However, we can prove that
the conditions in Proposition~\ref{prop:avg-op} and the hypotheses in
Theorems~\ref{thm:linear} and~\ref{thm:main} are related in two
important examples of Banach function spaces---the weighted and
variable Lebesgue spaces---and for singular integrals of convolution
type and fractional integral operators.

Before considering these spaces, we want to specify  the operators we are
interested in.   In the linear setting, we will consider singular
integrals  of the form
\[ Tf(x) = \text{p.v.}\int_{\nR^n} \frac{\Omega(y')}{|y|^n}f(x-y)\,dy, \]
where $x'=x/|x|^n$ and $\Omega $ is defined on $S^{n-1}$,
has mean  $0$, and is sufficiently smooth.  
Examples include the Riesz transforms $R_j$, which have kernels
$K_j(x) = \frac{x_j}{|x|^{n+1}}$.    For $0<\alpha<n$,  we will
consider the fractional
integral operator:  i.e., the convolution operator
\[ I_\alpha f(x) = \int_{\nR^n} \frac{f(y)}{|x-y|^{n-\alpha}}\,dy. \]
For more information on both kinds of operators,
see~\cite{grafakos08a,grafakos08b}.

In the bilinear setting, we consider singular integral operators of
the form
\[ T(f,g)(x) = \text{p.v.}\int_{\nR^{n}} \int_{\nR^{n}} 
\frac{\Omega((y_1,y_2)')}{|(y_1,y_2)|^n}f(x-y_1)g(x-y_2)\,dy_1dy_2, \]
where $\Omega$ is defined on $S^{2n-1}$, has mean $0$ and is
sufficiently smooth.  Examples include the multilinear Riesz
transforms.  For more on these operators,  see~\cite{MR1880324}.
We also consider the bilinear fractional integral operator, which is defined for $0<\alpha<2n$
by
\[ I_\alpha(f,g)(x) = \int_{\nR^{n}} \int_{\nR^{n}} 
\frac{f(y_1)g(y_2)}{(|x-y_1|+|x-y_2|)^{2n-\alpha}}\,dy_1dy_2. \]
For more on these operators see~\cite{MR1164632,MR2514845}.  

For brevity, in the following sections we will refer to linear and bilinear singular
integral operators whose kernels satisfy these hypotheses as regular
operators. 

\subsection*{Weighted Lebesgue spaces}

In this section we  apply Theorems~\ref{thm:linear}
and~\ref{thm:main} to the weighted Lebesgue spaces.  Given a weight $w$ (i.e.,
a non-negative, locally integrable function) and $1<p<\infty$, we
define the space $L^p(w)$ to be the set of all measurable functions $f$ such that
\[ \|f\|_{L^p(w)} = \left(\int_{\nR^n} |f(x)|^p
  w(x)\,dx\right)^{\frac{1}{p}} < \infty. \]
We say that a weight $w$ is in the Muckenhoupt class $A_p$ if for
every cube $Q$,
\[  \avgint_Q w(x)\,dx \left(\avgint_Q w(x)^{1-p'}\,dx\right)^{p-1}
\leq C.  \]
Then $L^p(w)$ is a Banach function space, and it is well known that
the associate space is $L^{p'}(w^{1-p'})$.  Further,  the
Hardy-Littlewood maximal operator is bounded on $L^p(w)$ if and only
if $w\in A_p$.  

For commutators, if $w\in A_p$ and
if $T$ is any Calder\'on-Zygmund singular integral operator (and not
just the class of singular integrals described above), and if
$b\in BMO$, then the commutator $[b,T] : L^p(w) \rightarrow L^p(w)$~\cite{perez97}.
Moreover, it is easy to see that the $A_p$ condition is equivalent to
\[ |Q|^{-1}\|\chi_Q\|_{L^p(w)}\|\chi_Q\|_{L^{p'}(w^{1-p'})} \leq C, \]
which is condition~\eqref{eqn:linear1}.   Therefore, we have proved
the following.

\begin{corollary} 
For $1<p<\infty$ and $w\in A_p$, given a regular singular
integral operator $T$ and a function $b$, if $[b,T]$ is bounded on
$L^p(w)$, then $b\in BMO$.
\end{corollary}

For $0<\alpha<n$ the corresponding weight condition
is the $A_{p,q}$ condition.  Given $1<p<\frac{n}{\alpha}$, define $q$
by $\frac{1}{p}-\frac{1}{q}=\frac{\alpha}{n}$.  We say $w\in A_{p,q}$
if for every cube $Q$,
\[ \left(\avgint_Q w(x)^q\,dx\right)^{\frac{1}{q}} \left(\avgint_Q
  w(x)^{-p'}\,dx\right)^{\frac{1}{p'}} \leq C. \]
We have that the fractional maximal operator satisfies $M_\alpha :
L^p(w^p)\rightarrow L^q(w^q)$ if and only if $w\in
A_{p,q}$~\cite{muckenhoupt-wheeden74}. 

For commutators of the fractional integral operator $I_\alpha$, 
if $w\in A_{p,q}$ and $b\in
BMO$,
$[b,I_\alpha] :  L^p(w^p) \rightarrow L^q(w^q)$ \cite{cruz-uribe-fiorenza03}.
The $A_{p,q}$ condition also implies \eqref{eqn:linear1}, though
unlike the case of $A_p$ weights, this is less obvious.  In this
case we have $X=L^p(w^p)$ and $Y=L^q(w^q)$, so $Y'=L^{q'}(w^{-q'})$,
and 
we can rewrite \eqref{eqn:linear1} as 
\[
\left( \avgint_Q w^{-q'}\,dx \right)^{\frac{1}{q'}}
\left( \avgint_Q w^{p}\,dx \right)^{\frac{1}{p}} \leq C. 
\]
Here we use the fact that 
since $\frac{1}{p}-\frac{1}{q}=\frac{\alpha}{n}$,
$\frac{1}{p}+\frac{1}{q'}=\frac{\alpha}{n}$.  Since $p<q$, $q'<p'$, so
if we apply H\"older's inequality twice we get that 
\[ \left( \avgint_Q w^{-q'}\,dx \right)^{\frac{1}{q'}}
\left( \avgint_Q w^{p}\,dx \right)^{\frac{1}{p}}
\leq \left( \avgint_Q w^{-p'}\,dx \right)^{\frac{1}{p'}}
\left( \avgint_Q w^{q}\,dx \right)^{\frac{1}{q}}
\leq C.\]

\begin{corollary}
  Given $0<\alpha<n$ and $1<p<\frac{n}{\alpha}$, define $q$ by
  $\frac{1}{p}-\frac{1}{q}=\frac{\alpha}{n}$.  Given $w\in A_{p,q}$
  and a function $b$, if the commutator
  $[b,I_\alpha] : L^p(w^p)\rightarrow L^q(w^q)$, then $b\in BMO$.
\end{corollary}

\medskip

We have similar results for bilinear operators, but they are much less
complete.   Given $1<p_1,\,p_2<\infty$, define the vector ``exponent''
$\vec{p}=(p_1,p_2,p)$, where
$\frac{1}{p}=\frac{1}{p_1}+\frac{1}{p_2}$.  Given $\vec{p}$ and
weights $w_1,\,w_2$, define the triple $\vec{w}=(w_1,w_2,w)$, where
$w=w_1^{\frac{p}{p_1}}w_2^{\frac{p}{p_2}}$.  Define the multilinear
analog of the Muckenhoupt $A_p$ weights as follows:  given $\vec{p}$,
we say that $\vec{w}\in A_{\vec{p}}$ if for every cube $Q$,
\[ \left(\avgint_Q w\,dx\right)^{\frac{1}{p}}
\left(\avgint_Q w_1^{1-p_1'}\,dx\right)^{\frac{1}{p_1'}}
\left(\avgint_Q w_2^{1-p_2'}\,dx\right)^{\frac{1}{p_2'}} \leq C. \]
These weights were introduced in~\cite{MR2483720}, where they showed
that the bilinear maximal
operator satisfies $M :L^{p_1}(w_1)\times L^{p_2}(w_2)\rightarrow
L^p(w)$ if and only if $\vec{w} \in A_{\vec{p}}$.   It is an immediate
consequence of H\"older's inequality that if $w_1\in A_{p_1}$ and
$w_2\in A_{p_2}$, then $\vec{w}\in A_{\vec{p}}$; however, this
condition is not necessary.

Given a bilinear Calder\'on-Zygmund singular integral
operator $T$ and $b\in BMO$,  we have that  if
$\vec{w}\in A_{\vec{p}}$, then for $i=1,2$,
$[b,T]_i : L^{p_1}(w_1)\times L^{p_2}(w_2)\rightarrow
L^p(w)$~\cite{MR2483720}.    In light of the results in the linear case, it
seems reasonable to conjecture that when $p>1$, the $A_{\vec{p}}$
condition implies~\eqref{eqn:main1}, which in this case can be written
as 
\begin{equation} \label{eqn:main-bilinear}
 \left(\avgint_Q w^{1-p'}\,dx\right)^{\frac{1}{p'}}
 \left(\avgint_Q w_1\,dx\right)^{\frac{1}{p_1}}
 \left(\avgint_Q w_2\,dx\right)^{\frac{1}{p_2}} \leq C. 
\end{equation}
(Here we use the fact that
$\frac{1}{p_1}+\frac{1}{p_2}+\frac{1}{p'}=1$.)  Written in this form,
this condition can be viewed formally as the bilinear analog of the fact that in the
linear case, if $w\in A_p$, then $w^{1-p'} \in A_{p'}$.  

However, we cannot prove this in general, or even in the special case
when $w_i\in A_{p_i}$, $i=1,2$.  We can prove that~\eqref{eqn:main1}
holds if we make the additional, stronger assumption that $w\in A_p$.
This holds, for instance, if we assume that
$w_1,\,w_2\in A_p\subset A_{p_i}$, $i=1,2$. (This inclusion holds since $A_q\subset A_r$ whenever $q<r$, and here
$p<p_i$.)  For then, since $w_i\in A_{p_i}$, by a multi-variable reverse H\"older
inequality recently proved in~\cite{DCU-KM} (and implicit
in~\cite[Theorem~2.6]{cruz-uribe-neugebauer95}), we have that
\[ \left(\avgint_Q w_1\,dx\right)^{\frac{p}{p_1}}
\left(\avgint_Q w_2\,dx\right)^{\frac{p}{p_2}} 
\lesssim \avgint_Q w\,dx. \]
But then, since $w\in A_p$,
\begin{multline*}
 \left(\avgint_Q w^{1-p'}\,dx\right)^{\frac{1}{p'}}
 \left(\avgint_Q w_1\,dx\right)^{\frac{1}{p_1}}
 \left(\avgint_Q w_2\,dx\right)^{\frac{1}{p_2}} \\
\approx \left(\avgint_Q w\,dx\right)^{-\frac{1}{p}}
 \left(\avgint_Q w_1\,dx\right)^{\frac{1}{p_1}}
 \left(\avgint_Q w_2\,dx\right)^{\frac{1}{p_2}}  \leq C.
\end{multline*}

\begin{corollary} \label{cor:wtd-sio-bilinear} Given $\vec{p}$ with
  $p>1$ and $\vec{w}\in A_{\vec{p}}$, suppose $w_i \in A_{p_i}$,
  $i=1,2$, and suppose $w\in A_p$.  If
  $T$ is a regular bilinear singular integral, and
  $b$ is function such that for $i=1,2$, $[b,T]_i : L^{p_1}(w_1)
  \times L^{p_2}(w_2) \rightarrow L^p(w)$, then $b\in BMO$.
\end{corollary}

\medskip

For bilinear fractional integrals, the corresponding weight class was
introduced in~\cite{MR2514845}.  With the notation as before, given
$0<\alpha<2n$ and 
$\vec{p}$, suppose that $\frac{1}{2}<p<\frac{n}{\alpha}$.  Define $q$ by 
  $\frac{1}{p}-\frac1q= \frac\alpha n$.  If we define the vector
  weight $\vec{w}=(w_1,w_2,w)$, where now $w=w_1w_2$, then
  $\vec{w}\in A_{\vec{p},q} $ if
\[\left(\avgint_Q w^{q}\,dx\right)^{\frac{1}{q}}  
\left(\avgint_Q w_1^{-p_1'}\,dx\right)^{\frac{1}{p_1'}}
\left(\avgint_Q w_2^{-p_2'}\right)^{\frac{1}{p_2'}} \leq C. \]
The bilinear fractional maximal operator satisfies $M_\alpha :
L^{p_1}(w^{p_1}) \times L^{p_2}(w_2) \rightarrow L^q(w^q)$ if and only
if $\vec{w} \in A_{\vec{p},q}$.  Similar to the $A_{\vec{p}}$ weights,
if $w_i \in A_{p_i,q_i}$, where $q_i>p_i$, $i=1,2$, and 
$\frac{1}{q_1}+\frac{1}{q_2}=\frac{1}{q}$, then $\vec{w} \in
A_{\vec{p},q}$.  

For the commutators of the bilinear fractional integral operator, 
if $w\in A_{\vec{p},q}$, then for $i=1,2$,
$[b,I_\alpha]_i: L^{p_1}(w_1^{p_1})\times L^{p_2}(w_2^{p_2})\to
L^q(w^{q})$ \cite{CW,CX}.  
As we did for singular integrals, we conjecture that the
$A_{\vec{p},q}$ condition implies~\eqref{eqn:main1}, which in this
case can be written as 
\begin{equation} \label{eqn:bilinear-frac}
 \left(\avgint_Q w^{-q'}\,dx\right)^{\frac{1}{q'}}
 \left(\avgint_Q w_1^{p_1}\,dx\right)^{\frac{1}{p_1}}
 \left(\avgint_Q w_2^{p_2}\,dx\right)^{\frac{1}{p_2}} \leq C. 
\end{equation}
(Here we use the fact that
$\frac{1}{p_1}+\frac{1}{p_2}+\frac{1}{q'}=1+\frac{\alpha}{n}$.)  

Arguing as in the case of bilinear fractional singular integrals, we
can prove this if we assume that $w_i \in A_{p_i,q_i}$ and $w^q \in A_q$:  i.e., that 
\[  \left(\avgint_Q w^{q}\,dx\right)^{\frac{1}{q}} \left(\avgint_Q
    w^{-q'}\,dx\right)^{\frac{1}{q'}} \leq C.   \]
For then, again by
the bilinear reverse H\"older inequality (since $w\in
A_{p_i,q_i}$,  $w\in A_\infty$, and this is sufficient for this
inequality to
hold) and H\"older's
inequality (since $q_i>p_i$),
\[ \left(\avgint_Q w^{q}\,dx\right)^{\frac{1}{q}}
\geq 
\left(\avgint_Q w_1^{q_1}\,dx\right)^{\frac{1}{q_1}}
 \left(\avgint_Q w_2^{q_2}\,dx\right)^{\frac{1}{q_2}} 
\geq
\left(\avgint_Q w_1^{p_1}\,dx\right)^{\frac{1}{p_1}}
 \left(\avgint_Q w_2^{p_2}\,dx\right)^{\frac{1}{p_2}};\]
Inequality~\eqref{eqn:bilinear-frac} follows at once.

We can eliminate the assumption that $w^q\in A_q$ if we
restrict the range of $\alpha$ to $n \leq \alpha<2n$. 
Suppose $w_i\in A_{p_i,q_i}$, $i=1,2$; then 
 we have that $\frac{1}{q} = \frac{1}{p}-\frac{\alpha}{n}
\leq \frac{1-p}{p} < 1$ since $p>\frac{1}{2}$.  Moreover, we have that
$\frac1{q'}=\frac1{p'_1}+\frac1{
  p_2'}+\left(\frac\alpha n-1\right)\geq \frac1{p'_1}+\frac1{ p_2'}$.
Therefore, we can apply H\"older's inequality three times to get that 
the left-hand side of \eqref{eqn:bilinear-frac} is bounded by
\[  \left(\avgint_Q w_1^{-p_1'}\,dx\right)^{\frac{1}{p_1'}}
 \left(\avgint_Q w_2^{-p_2'}\,dx\right)^{\frac{1}{p_2'}} 
\left(\avgint_Q w_1^{q_1}\,dx\right)^{\frac{1}{q_1}}
 \left(\avgint_Q w_2^{q_2}\,dx\right)^{\frac{1}{q_2}}\leq C.  \]

\begin{corollary} \label{cor:wtd-frac-bilinear} Given $0<\alpha<2n $,
  $\vec{p}$ such that $p<\frac{n}{\alpha}$, and $\vec{w}$ such that
  $w_i \in A_{p_i,q_i}$ with $q_i>p_i$, suppose either $w^q\in A_q$ or
  $n\leq \alpha<2n$.   If $b$ is a function such that for $i=1,2$,
  $[b,I_\alpha]_i : L^{p_1}(w_1^{p_1}) \times L^{p_2}(w_2^{p_2})
  \rightarrow L^q(w^q)$, then $b\in BMO$.
\end{corollary}

\begin{remark}
In Corollaries~\ref{cor:wtd-sio-bilinear}
and~\ref{cor:wtd-frac-bilinear}, we can interpret the
hypotheses $w\in L^p(w)$ and $w^q\in L^q(w^q)$ as assuming the maximal operator is bounded on the target
space $L^p(w)$ or $L^q(w^q)$.  This should be compared to the
assumptions in Corollaries~\ref{cor:var-sio-bilinear}
and~\ref{cor:var-frac-bilinear} below.
\end{remark}


\subsection*{Variable Lebesgue spaces}
The variable Lebesgue spaces are a generalization of the classical
$L^p$ spaces.  
Given a measurable function $\pp : \nR^n \rightarrow
[1,\infty)$, we define $L^\pp$ to be the collection of all measurable
functions such that
\[ \|f\|_{L^\pp} =\inf\left\{ \lambda > 0 :
\int_{\nR^n} \bigg(\frac{|f(x)|}{\lambda}\bigg)^{p(x)}\,dx \leq 1
\right\}. \]
With this norm $L^\pp$ is a Banach function space; the associate space
is $L^\cpp$, where we define $\cpp$ pointwise by
$\frac{1}{p(x)}+\frac{1}{p'(x)}=1$.     For brevity, we define
\[ p_- = \essinf_{x\in \nR^n} p(x), \qquad p_+ = \esssup_{x\in \nR^n}
p(x). \]
For complete information on these spaces, see~\cite{cruz-fiorenza-book}.  

The boundedness of the maximal operator depends (in a very subtle way)
on the regularity of the exponent function $\pp$.  
A sufficient condition for $M : L^\pp \rightarrow L^\pp$ is that $p_- >1$ and $\pp$ satisfies the
log-H\"older continuity conditions locally and at infinity:  
\[ \left|\frac{1}{p(x)}-\frac{1}{p(y)}\right| \leq
\frac{C_0}{-\log(|x-y|)}, \quad |x-y|\leq \frac{1}{2}, \]
and there exists a constant $1\leq p_\infty \leq \infty$ such that
\[ \left|\frac{1}{p(x)} - \frac{1}{p_\infty}\right|
\leq \frac{C_\infty}{\log(e+|x|)}. \]
By Proposition~\ref{prop:avg-op}, if $M : L^\pp
\rightarrow L^\pp$, then  for every cube~$Q$,
\begin{equation} \label{eqn:K0}
\|\chi_Q\|_{L^\pp}\|\chi_Q\|_{L^\cpp} \leq C|Q|, 
\end{equation}
which is \eqref{eqn:linear1}.  However, we have a stronger result.
Suppose $\pp$ is such that
the maximal operator is bounded on $L^\pp$.  Given a cube $Q$, if we define the
exponents $p_Q$ and $p'_Q$ by
\[ \frac{1}{p_Q} = \avgint_Q \frac{1}{p(x)}\,dx, \qquad 
\frac{1}{p'_Q} = \avgint_Q \frac{1}{p'(x)}\,dx, \]
then 
\begin{equation} \label{eqn:strong-K0}
\|\chi_Q\|_\pp \approx |Q|^{\frac{1}{p_Q}} \quad \text{ and } \quad \|\chi_Q\|_\cpp
\approx |Q|^{\frac{1}{p'_Q}},
\end{equation}
and the implicit constants are
independent of $Q$~\cite[Proposition~4.66]{cruz-fiorenza-book}.

Let $T$ be a Calder\'on-Zygmund singular integral operator.  If $\pp$ is such that
$1<p_-\leq p_+<\infty$ and the maximal operator is bounded on $L^\pp$,
then for all $b\in BMO$,  $[b,T] :  L^\pp \rightarrow
L^\pp$~\cite[Corollary~2.10]{MR2210118}.

\begin{corollary}
Let $\pp$ be an exponent function such that $1<p_-\leq p_+<\infty$ and
the maximal operator is bounded on $L^\pp$.  If $T$ is a regular
singular integral and $b$ is a function such that $[b,T]$ is bounded
on $L^\pp$, then $b\in BMO$.  
\end{corollary}

Given $0<\alpha<n$ and $\pp$  such that $1<p_-\leq p_+ <
\frac{n}{\alpha}$, define $\qq$ pointwise by
$\frac{1}{p(x)}-\frac{1}{q(x)}=\frac{\alpha}{n}$.    If there
exists $q_0>\frac{n}{n-\alpha}$ such that the maximal operator is bounded
on $L^{(\qq/q_0)'}$, then for all $b\in BMO$, $[b,I_\alpha] : L^\pp
\rightarrow L^\qq$.  (This result does not appear explicitly in the
literature, but it is a straightforward application of known results.
For instance, it follows by extrapolation, arguing as in
\cite[Theorem~5.46]{cruz-fiorenza-book} but using the weighted norm
inequalities for commutators
from~\cite[Theorem~1.6]{cruz-uribe-fiorenza03}.) 

If the maximal operator is bounded on $L^{(\qq/q_0)'}$, then
by~\cite[Corollary~4.64]{cruz-fiorenza-book}, it is also bounded on
$L^{\qq/q_0}$ and so on $L^\qq$ \cite[Theorem~4.37]{cruz-fiorenza-book}.   If we let $\theta=\frac{1}{q_0}$, then we can write
\[ \frac{1}{p(x)} = \frac{1}{q(x)}+\frac{\alpha}{n}
= \frac{\theta}{q(x)/q_0} +
\frac{1-\theta}{(1-\theta)\frac{n}{\alpha}}. \]
By our assumption on $q_0$, $r=(1-\theta)\frac{n}{\alpha}>1$, and so
the maximal operator is bounded on $L^r$.  Hence, 
by interpolation (see~\cite[Theorem~3.38]{cruz-fiorenza-book}) the
maximal operator is bounded on $L^\pp$.  Therefore, by
\eqref{eqn:strong-K0},

\[ |Q|^{-\frac{\alpha}{n}}\|\chi_Q\|_{L^\cqq}\|\chi_Q\|_{L^\pp}
\lesssim
|Q|^{-\frac{\alpha}{n}} |Q|^{\frac{1}{q'_Q}}|Q|^{\frac{1}{p_Q}}
\lesssim 
|Q|. \]
So again, \eqref{eqn:linear1} holds.

\begin{corollary}
Given $0<\alpha<n$ and $\pp$  such that $1<p_-\leq p_+ <
\frac{n}{\alpha}$, define $\qq$ pointwise by
$\frac{1}{p(x)}-\frac{1}{q(x)}=\frac{\alpha}{n}$.    Suppose there
exists $q_0>\frac{n}{n-\alpha}$ such that the maximal operator is bounded
on $L^{(\qq/q_0)'}$.  If $b$ is such that $[b,I_\alpha] : L^\pp
\rightarrow L^\qq$, then $b\in BMO$.
\end{corollary}

\medskip

We have similar results for bilinear operators.  Suppose $p_1(\cdot)$,
$p_2(\cdot)$ are such that $1<(p_i)_-\leq (p_i)_+ < \infty$, $i=1,2$,
and define $\pp$ pointwise by
$\frac{1}{p(x)} = \frac{1}{p_1(x)}+\frac{1}{p_2(x)}$.  Note that
$p_->\frac{1}{2}$.  If we further assume that the (linear) maximal
operator is bounded on $L^{p_i(\cdot)}$, $i=1,2$, then given any
bilinear Calder\'on-Zygmund singular integral operator $T$, we have
that 
$T :  L^{p_1(\cdot)}\times L^{p_2(\cdot)} \rightarrow
L^\pp$~\cite[Corollary~4.1]{CruzUribe:2016wv}.  Moreover, given any
$b\in BMO$,
$[b,T]_i : L^{p_1(\cdot)}\times L^{p_2(\cdot)} \rightarrow L^\pp$.  This
is not proved explicitly in the literature, but the proof is
the same as for bilinear singular integrals, using bilinear extrapolation and
using the weighted inequalities for bilinear commutators
in~\cite{MR2483720}.  

With the same assumptions we also
have that $M : L^{p_1(\cdot)}\times L^{p_2(\cdot)} \rightarrow L^\pp$:
by the generalized H\"older's
inequality~\cite[Corollary~2.28]{cruz-fiorenza-book},
\[ \|M(f,g)\|_\pp \leq \|Mf \cdot Mg\|_\pp 
\leq \|Mf\|_{p_1(\cdot)}\|Mg\|_{p_2(\cdot)}. \]
Note that we pass from the bilinear to the linear maximal operator.
Also, note that the generalized H\"older's inequality is only proved
in~\cite{cruz-fiorenza-book} assuming $p_-\geq 1$, but for $p_->0$ it
follows by a rescaling argument:
cf.~\cite[Lemma~2.7]{CruzUribe:2014ux}. 

If we further assume $p_->1$ and the maximal operator is bounded on
$L^\pp$, then by~\eqref{eqn:strong-K0} we have that 
\[ \|\chi_Q\|_\cpp \|\chi_Q\|_{p_1(\cdot)}\|\chi_Q\|_{p_2(\cdot)}
\lesssim
|Q|^{\frac{1}{p'_Q}+\frac{1}{(p_1)_Q}+\frac{1}{(p_2)_Q}} \leq C.  \]

\begin{corollary} \label{cor:var-sio-bilinear}
Suppose $p_1(\cdot)$,
$p_2(\cdot)$ are such that $1<(p_i)_-\leq (p_i)_+ < \infty$, $i=1,2$, and
we define $\pp$ pointwise by $\frac{1}{p(x)} =
\frac{1}{p_1(x)}+\frac{1}{p_2(x)}$.   Suppose further that $p_->1$ and
the maximal
operator is bounded on $L^\pp$ and $L^{p_i(\cdot)}$, $i=1,2$.  If $T$
is a regular bilinear singular integral and $b$ is such that 
$[b,T]_i :  L^{p_1(\cdot)}\times L^{p_2(\cdot)}
\rightarrow L^\pp$, then $b\in BMO$.
\end{corollary}

To prove the analogous result for the bilinear fractional integral
operator, fix $0<\alpha<2n$.  Suppose $p_1(\cdot)$, $p_2(\cdot)$ are
such that $1<(p_i)_-\leq (p_i)_+ < \infty$, $i=1,2$, and again define
$\pp$ pointwise by
$\frac{1}{p(x)} = \frac{1}{p_1(x)}+\frac{1}{p_2(x)}$.  Define $\qq$ by
$\frac{1}{p(x)}-\frac{1}{q(x)}=\frac{\alpha}{n}$.  Fix
$0<\alpha_1,\alpha_2<n$ such that $\alpha_1+\alpha_2=\alpha$, and define
$q_i(\cdot)$, $i=1,2$, by
$\frac{1}{p_i(x)}-\frac{1}{q_i(x)}=\frac{\alpha_i}{n}$.  If we assume
that the maximal operator is bounded on $L^\qq$, and that there exist
$q_i>\frac{n}{n-\alpha_i}$. $i=1,2$, such that the maximal operator is
bounded on $L^{(q_i(\cdot)/q_i)'}$, then given any $b\in BMO$,
$[b,I_\alpha]_i : L^{p_1(\cdot)}\times L^{p_2(\cdot)} \rightarrow
L^\qq$.  As in the linear case, this result has not been explicitly
proved in the literature, but follows from known results.  For all
$w\in A_\infty$ and any $0<p<\infty$,
$\|I_\alpha(f,g)\|_{L^p(w)}\lesssim
\|M_\alpha(f,g)\|_{L^p(w)}$~\cite[Theorem~3.1]{MR2514845}.  Since the
maximal operator is bounded on $L^\qq$, by extrapolation,
$\|I_\alpha(f,g)\|_\qq \lesssim
\|M_\alpha(f,g)\|_\qq$~\cite[Theorem~5.24]{cruz-fiorenza-book}.  By
the generalized H\"older's inequality,
\[ \|M_\alpha(f,g)\|_\qq 
\leq \|M_{\alpha_1}f \cdot M_{\alpha_2}g\|_\qq
\leq \|M_{\alpha_1}f\|_{q_1(\cdot)} \|M_{\alpha_2}g\|_{q_2(\cdot)} 
\leq \|f\|_{p_1(\cdot)} \|g\|_{p_2(\cdot)} ; \]
the last inequality follows from our assumptions on $q_i(\cdot)$
and~\cite[Remark~5.51]{cruz-fiorenza-book}.

In this case~\eqref{eqn:main1} becomes
\[ |Q|^{-\frac{\alpha}{n}}\|\chi_Q\|_\cqq \|\chi_Q\|_{p_1(\cdot)}
\|\chi_Q\|_{p_2(\cdot)} \leq C|Q|.  \]
If we make the same assumptions on the exponents as used to prove the
inequalities for the commutators, then arguing as we did above
using~\eqref{eqn:strong-K0}, we get this inequality.

\begin{corollary} \label{cor:var-frac-bilinear}
 Fix $0<\alpha<2n$.  Suppose $p_1(\cdot)$,
$p_2(\cdot)$ are such that $1<(p_i)_-\leq (p_i)_+ < \infty$, $i=1,2$, and
again define $\pp$ pointwise by $\frac{1}{p(x)} =
\frac{1}{p_1(x)}+\frac{1}{p_2(x)}$.  Define $\qq$ by
$\frac{1}{p(x)}-\frac{1}{q(x)}=\frac{\alpha}{n}$.  Fix
$0<\alpha_1,\alpha_2<n$ such that $\alpha_1+\alpha_2=\alpha$, and define
$q_i(\cdot)$, $i=1,2$, by
$\frac{1}{p_i(x)}-\frac{1}{q_i(x)}=\frac{\alpha_i}{n}$.  Suppose
that the maximal operator is bounded on $L^\qq$, and that there exist
$q_i>\frac{n}{n-\alpha_i}$. $i=1,2$, such that the maximal operator is
bounded on $L^{(q_i(\cdot)/q_i)'}$.  Given any $b$, if $[b,I_\alpha]_i : L^{p_1(\cdot)}\times L^{p_2(\cdot)}
\rightarrow L^\qq$, then $b\in BMO$.
\end{corollary}

\medskip

\begin{remark}
  Sufficient conditions for the boundedness of commutators of singular and fractional
  integrals on Orlicz spaces are known or can be readily proved using
  extrapolation:
  see~\cite{curbera-garcia-cuerva-martell-perez06,kokilashvili-krbec91}.  
  Similarly, such results can be proved in the weighted variable Lebesgue spaces and
  generalized Orlicz spaces (also known as Nakano spaces or
  Musielak-Orlicz spaces) using extrapolation:
  see~\cite{DCU-PH,CruzUribe:2015km} for definitions and the corresponding extrapolation results.  Results
  similar to those for the variable Lebesgue spaces above can be
  deduced in these settings---we leave the precise statements and
  proofs to the interested reader.
\end{remark}

\section{Proof of Theorem~\ref{thm:main}}
\label{section:proof}

In this section we prove Theorem~\ref{thm:main}.  As we noted in the
Introduction, the proof of Theorem~\ref{thm:linear} is gotten by a
simple adaptation of this proof.   Our proof actually closely follows the argument due
to the first author in~\cite{LC}, which was in turn following the
techniques of Janson in \cite{SJ}.     Here we will concentrate on the
parts of the proof which changes because we are working in the setting
of Banach function spaces, and we refer the reader to~\cite{LC} for
further details.

\begin{proof}
  We will assume without loss of generality that $[b,T]_1$ is bounded;
  the proof for the other commutator is identical.

Let $B=B((y_0,z_0),\delta\sqrt{2n})\subset\nR^{2n}$ be a ball upon
  which $\frac1K$ has an absolutely convergent Fourier series.  By the
  homogeneity of $K$, we may assume without loss of generality that
  $2\sqrt n<|(y_0,z_0)|<4\sqrt n$ and $\delta<1$. These conditions
  guarantee that $\overline B \cap \{0\}=\emptyset$, avoiding any
  potential singularity of $K$; {this will be important below as it
    will let us use the integral representation of $[b,T]_1$.}

Write the Fourier series of $\frac{1}{K}$ as
\[\frac1{K(y,z)}=\sum_ja_je^{i\nu_j\cdot(y,z)}=\sum_ja_je^{i(\nu_j^1,\nu_j^2)\cdot(y,z)};\] 
note that the individual vectors
$\nu_j=(\nu_j^1,\nu_j^2)\in\nR^n\times\nR^n$ do not play any
significant role in the proof, and we introduce them simply to be precise.

Let $y_1=\delta^{-1}y_0$, $z_1=\delta^{-1}z_0$; then by homogeneity we
have that for all $(y,z)\in B((y_1,z_1),\sqrt{2n})$, 
\[ \frac1{K(y,z)}=\frac{\delta^{-2n+\alpha}}{K(\delta y,\delta z)}=\delta^{-2n+\alpha}\sum_ja_je^{i\delta\nu_j\cdot(y,z)}.\]
Let $Q=Q(x_0,r)$ be an arbitrary cube in $\nR^n$, and set $\tilde
y=x_0+ry_1,\ \tilde z=x_0+rz_1$.  Define $Q'=Q(\tilde y,r)$ and
$Q''=Q(\tilde z,r)$.

It follows from the size conditions on $y_0$ and $z_0$ that $Q$ and
either $Q'$ or $Q''$ are disjoint.  To see that this is the case, note that the
minimum size condition on $(y_0,z_0)$ implies that
$\max\{|y_0|,|z_0|\}\geq\sqrt{2n}$; without loss of generality,
suppose that it is $|y_0|$. This in turn implies that the distance
between $x_0$ and $\tilde y$ is greater than $r\sqrt{2n}$. Since $Q$
and $Q'$ each have side-length $r$, the distance of their centers from
one another guarantees that they must be disjoint. If $|z_0|$ is
larger, then we get the same conclusion for $Q''$. 

As a consequence, $Q\cap Q'\cap Q''=\emptyset$, which allows us to use
the kernel representation of $[b,T]_1$ for
$(x,y,z)\in Q\times Q'\times Q''$.  Additionally, for
$(x,y,z)\in Q\times Q'\times Q''$,
$\left(\frac{x-y}r,\frac{x-z}r\right)\in B((y_1,z_1),\sqrt{2n})$,
which in turn means that $(x-y,x-z)$ is bounded away from the
singularity of $K$, and so the integral representation of $[b,T]_1$ can
be used freely. Further details of these calculations can be found in
\cite{LC}. 

 It also follows from the maximum size
condition on $(y_0,z_0)$ that
\begin{equation} \label{eqn:size-cubes}
Q',\; Q''\subset B\left(x_0,\left(\frac{r\sqrt n}{2}\left(1+\frac8\delta\right)\right)\right)\subset\sqrt
n\left(1+\frac8\delta\right)Q.
\end{equation}
To see this, note that the maximum size of $y_0$ or $z_0$ is $4\sqrt
n$ which implies that the maximum distance from $x_0$ to $\tilde y$ or
$\tilde z$  is $\frac{4r\sqrt n}{\delta}$.  The final containment
in $B\left(x_0,\left(\frac{r\sqrt
      n}{2}\left(1+\frac8\delta\right)\right)\right)$ follows from this.

\medskip

We can now estimate as follows.  Fix $Q$ and  let
$\sigma(x)=\sgn(b(x)-b_Q)$; then
\begin{align*}
\int_Q|b(x)&-b_{Q'}|\,dx\\
&=\int_Q(b(x)-b_{Q'})\sigma(x)\,dx\\
&=\frac1{|Q''|}\frac1{|Q'|}\int_Q\int_{Q'}\int_{Q''}(b(x)-b(y))\sigma(x)\,dzdydx\\
&=r^{-2n}\int_{\nR^n}\int_{\nR^n}\int_{\nR^n}(b(x)-b(y))\frac{r^{2n-\alpha}K(x-y,x-z)}{K\left(\frac{x-y}r,\frac{x-z}r\right)}\\
&\quad\quad\times\sigma(x)\chi_{Q}(x)\chi_{Q'}(y)\chi_{Q''}(z)\,dzdydx\\
&=\delta^{-2n+\alpha}r^{-\alpha}\int_{\nR^n}\int_{\nR^n}\int_{\nR^n} 
(b(x)-b(y))K(x-y,x-z)\sum_ja_je^{i\frac\delta r\nu_j\cdot(x-y,x-z)}\\
&\quad\quad\times\sigma(x)\chi_{Q}(x)\chi_{Q'}(y)\chi_{Q''}(z)\,dzdydx
\end{align*}

Define the functions
\begin{align*}
f_j(y)&=e^{-i\frac\delta r\nu^1_j\cdot y}\chi_{Q'}(y)\\
g_j(z)&=e^{-i\frac\delta r\nu^2_j\cdot z}\chi_{Q''}(z)\\
h_j(x)&=e^{i\frac\delta r\nu_j\cdot(x,x)}\sigma(x)\chi_{Q}(x).
\end{align*}
Note that the norm of each of these functions in any Banach function
space will be the same as the norm of its support.   We can now
continue the above estimate:
\begin{align*}
\int_Q|b(x)-&b_{Q'}|~dx\\
&=\delta^{-2n+\alpha}r^{-\alpha}\sum_ja_j
\int_{\nR^n} h_j(x)\int_{\nR^n}\int_{\nR^n} (b(x)-b(y))\\
&\quad\quad\times K(x-y,x-z)f_j(y)g_j(z)~dzdydx\\
&=\delta^{-2n+\alpha}|Q|^{-\frac{\alpha}{n}}\sum_ja_j
\int_{\nR^n} h_j(x)[b,T]_1(f_j,g_j)(x)~dx\\
&\leq\delta^{-2n+\alpha}|Q|^{-\frac{\alpha}{n}}\sum_j|a_j|
\int_{\nR^n} |h_j(x)||[b,T]_1(f_j,g_j)(x)|~dx\\
&\leq\delta^{-2n+\alpha}|Q|^{-\frac{\alpha}{n}}\sum_j|a_j|\|h_j\|_{Y'}\|[b,T]_1(f_j,g_j)\|_{Y}\\
&\leq\delta^{-2n+\alpha}|Q|^{-\frac{\alpha}{n}}\sum_j|a_j|\|h_j\|_{Y'}\|[b,T]_1\|_{X_1\times X_2\to Y}\|f_j\|_{X_1}\|g_j\|_{X_2}\\
&=\delta^{-2n+\alpha}\|[b,T]\|_{X_1\times X_2\to Y}\sum_j|a_j|\|\chi_Q\|_{Y'}\|\chi_{Q'}\|_{X_1}\|\chi_{Q''}\|_{X_2}|Q|^{-\frac{\alpha}{n}}.\\
\end{align*}

Let $P=2\sqrt n(1+\frac8\delta)Q$.  By
inequality~\eqref{eqn:size-cubes} we have that  $Q,\,Q',\,Q"\subset
P$, and $|P|\approx |Q|$.  Therefore, by~\eqref{eqn:main1},
\begin{align*}
\int_Q|b(x)-&b_{Q'}|~dx\\
&\lesssim\|[b,T]\|_{X_1\times X_2\to Y}\sum_j|a_j|\|\chi_P\|_{Y'}\|\chi_{P}\|_{X_1}\|\chi_{P}\|_{X_2}|P|^{-\frac{\alpha}{n}}\\
&\lesssim|P|\|[b,T]_1\|_{X_1\times X_2\to Y}\sum_j|a_j|, \\
& \lesssim |Q|.
\end{align*}
Since this is true for every cube $Q$, $b\in BMO$ and the proof is complete.
\end{proof}

\bibliographystyle{plain}
\bibliography{bmo-necessary}

\end{document}